\documentclass{article}
\usepackage{graphicx}
\usepackage{amsthm, amsmath, amssymb}
\newtheorem{thm}{Theorem}[section]

\newtheorem{rem}[thm]{Remark}
\numberwithin{equation}{section}

\begin{document}



\title{Bounds on discriminants with one class per genus}%
\author{Kimberly Hopkins\\UT Austin, Department of Mathematics C1200\\Austin, TX 78712\\
\date{January 23, 2004}
\texttt{khopkins@math.utexas.edu }}%
 
\maketitle

\begin{abstract}
We assume a bound on the class number proposed by Conrey and Iwaniec in \cite{CI} to show there are no
discriminants with one class per genus beyond $d_{66} \approx 1.9 \cdot 10^{130}$.
\end{abstract}
\maketitle
\section{Introduction}

Let $-d<0$ be the discriminant of the complex quadratic field $\mathbb{Q}(\sqrt{-d})$ with ideal class group
$\mathcal{C}(-d)$ and class number $h(-d)=\sharp\:\mathcal{C}(-d)$. Recall the corresponding reduced quadratic
forms are subdivided into precisely $2^{g-1}$ genera where $g$ is the number of primary discriminants dividing
$d$.

The discriminants with one class per genus are of particular interest. The class group for such a discriminant
is completely determined
\begin{equation*}
    \mathcal{C}(-d) \cong (\mathbb{Z}/2 \mathbb{Z})^{g-1}
\end{equation*}
 and the class number is  $h(-d) = 2^{g-1}$.

Discriminants with one class per genus which are congruent to $0$ modulo $4$ are Euler's Idoneal (or convenient)
numbers which he used in the computation of large primes. There are $65$ known examples of such discriminants
\cite[p.59-63]{Cox}.

 In \cite[p.118]{Wein} Weinberger showed there to be at most one more. Assuming GRH for $L(s,\chi)$, Weinberger also
 proved there are
  no discriminants with one class per genus besides the $65$ known examples. (A $66$th would thus contradict GRH.)
   He did this using the bound
\begin{equation}\label{tatuzawa}
    L(1, \chi)> \frac{.655}{e}\frac{1}{\log d}
\end{equation}
given by T. Tatuzawa in \cite{Tat}.

In \cite[p.263]{CI} B. Conrey and H. Iwaniec assume a weaker condition on the \emph{spacing} between zeros on the
critical line to obtain other bounds on the class number, such as
\begin{equation}\label{conreyiwaniec}
    L(1,\chi) \geq \frac{1}{(\log d)^{18}}.
\end{equation}

This bound is not as sharp as Tatuzawa's. However, no assumptions regarding the GRH are required. On the
contrary, zeros off the critical line are permitted so long as the zeros on the line satisfy the proper spacing
criteria. In this note we will use Conrey and Iwaniec's bound to prove a result analogous to Weinberger's
concerning discriminants with one class per genus.

\section{Results of Conrey-Iwaniec}
Let $L(s,\psi) =\sum_{\mathfrak{a}} \psi (\mathfrak{a})(N\mathfrak{a})^{-s}$ be the L-function associated with the
class group character $\psi \in \hat{\mathcal{C}}(-d)$, where $\mathfrak{a}$ runs over the non-zero integral
ideals of $\mathbb{Q}(\sqrt{-d})$.   For any zero $\rho = \frac12 + i \gamma$ of $L(s,\psi)$ on the critical line,
let $\rho ' = \frac12+ i\gamma '$ denote the zero on the critical line closest to $\rho$. One can show that the
number of zeros of $L(s, \psi)$ in the rectangle $s=\sigma + i t$ with $0\leq \sigma \leq 1$, $0< t\leq T$
satisfies
\begin{equation*}
    N(T, \psi) = \frac{T}{\pi}\log \frac{T\sqrt{d}}{2\pi} - \frac{T}{\pi} + O(\log dT).
\end{equation*}
Assuming GRH for $L(s,\psi)$, this implies that the average gap between consecutive zeros $\rho$ and $\rho'$ is
about $\pi/\log\gamma$ \cite[p.263]{CI}.

Conrey and Iwaniec prove that if the gap is somewhat smaller than the average gap for sufficiently many pairs of
zeros (no Riemann hypothesis is required), then $h(-d) \gg \sqrt{d}(\log d)^{-A}$ for some constant $A \geq 0$.

 For $0< \alpha \leq 1$ and $T\geq 0$, define
\begin{equation}\label{pairs}
    D(\alpha, T) := \sharp \{\rho : 2\leq \gamma \leq T, |\gamma - \gamma '|\leq \frac{\pi (1-\alpha)}{\log
    \gamma}\}.
\end{equation}
 Then $D(\alpha, T) $ measures the number of consecutive zeros with smaller than average gaps. Note they do not
 count zeros off the critical line, but do allow them to exist.

Let $L(s,\chi) = \sum_{n=1}^{\infty} \chi(n)n^{-s}$ be the L-series for the real character $\chi$ of conductor
$d$. In Theorem $1.1$ of \cite{CI}, Conrey and Iwaniec prove
\begin{thm}\label{CI Theorem}
Let $A\geq 0$ and $\log T \geq (\log d)^{A+6}$. Suppose
\begin{equation}\label{conreyresult}
    D(\alpha, T) \geq \frac{c T (\log T)}{\alpha (\log d)^A}
\end{equation}
for some $0<\alpha \leq 1$, where $c$ is a large absolute constant. Then
\begin{equation}\label{conrey result2}
    L(1,\chi) \geq (\log T)^{-2}(\log d)^{-2A -6}.
\end{equation}
\end{thm}
This theorem gives a remarkable connection between clustering of zeros on the critical line and the class number.

\section{Theorem}
Analogous to Weinberger's assumption of GRH, we assume the weaker condition \eqref{pairs} to get a bound on
$h(-d)$. For a fixed $-d<0$ as in the previous section, we choose $A=0$ and $\log T=  (\log d)^6$ which under the
assumption of $(2.2)$ and Theorem \ref{CI Theorem} implies
 \begin{equation}\label{bestbound}
    L(1,\chi) \geq \frac{1}{(\log d)^{18}}.
\end{equation}
This is the best possible bound that can be obtained from Theorem \ref{CI Theorem}. With this in mind we
prove the following theorem:

\begin{thm}\label{My Theorem}
 Let $d_g$ denote the product of the first $g$ primes. Suppose $L(1, \chi) \geq (\log d)^{-18}$ for a fundamental
 discriminant $-d$. If $-d$ has one class
per genus then $d\leq d_{66}\approx 1.9 \cdot 10^{130}$.\\
 \end{thm}
 \begin{rem}
 Assuming all fundamental discriminants satisfy this bound
implies there are no discriminants with one class per genus beyond $d_{66}$.
 \end{rem}
 \begin{proof}
 By way of contradiction,
assume $d>d_{66}$ with one class per genus. \\

From the analytic class number formula and our hypothesis, we have \begin{equation*}
    h(-d)\geq \frac{\sqrt{d}}{\pi (\log d)^{18}}.
\end{equation*}
The derivative of $f(x):=\frac{\sqrt{x}}{\pi (\log x)^{18}}$ is $f'(x) = \frac{\log x - 36}{2\pi \sqrt{x}(\log
x)^{19}}$ which is zero only at $x=e^{36}$ and positive to its right. Hence $f(x)$ is an increasing function
for $x\geq e^{36}$.\\
We get
\begin{eqnarray*}
  2^{g-1} &=& h(-d) \\
   &\geq&  \frac{\sqrt{d}}{\pi (\log d)^{18}}\\
   &>& \frac{\sqrt{d_{66}}}{\pi (\log d_{66})^{18}} \\
   &>& 1.1 \cdot 10^{20}.
\end{eqnarray*}
This implies $g$ must be at least $68$ since $2^{66}<1.1 \cdot 10^{20} < 2^{67}$ and $2^{g-1}$ is an increasing
function in $g$.

In particular we may assume $g>66$. This gives us the inequality
\begin{equation}\label{inequality}
    d \geq d_g > d_{66}\cdot 331^{g-66}
\end{equation}
where  $331$ is the $67$th prime.

Using this inequality we have
\begin{eqnarray*}
  2^{g-1}=h(-d) &\geq& \frac{\sqrt{d}}{\pi (\log d)^{18}}\\
  &\geq& \frac{\sqrt{d_g}}{\pi (\log d_g)^{18}} \\
   &>& \frac{\sqrt{d_{66}}\cdot 331^{g-66/2}}{\pi[ \log d_{66}+(g-66)\log 331]^{18}}.
\end{eqnarray*}
Multiplying through by $2^{66-g}$ gives us
\begin{equation}\label{star}
    2^{65} >\frac{\sqrt{d_{66}}\cdot 9^{g-66}}{\pi[ \log d_{66}+(g-66)\log 331]^{18}}.
\end{equation}
Again letting
\begin{equation*}
r(x):=\frac{9^x}{[ \log d_{66}+x\log 331]^{18}}
\end{equation*}
 we see

 \begin{equation*}
    r'(x) = \frac{554.7\cdot 9^x + x \cdot 9^x \log 9 \log 331}{[300 + x\log 331]^{19}}
\end{equation*}
 is positive
for $x\geq 0$. Hence the right hand side of \eqref{star} is increasing for $g>66$ and so we can substitute
$g=67$ into \eqref{star}  to obtain
\begin{eqnarray*}
  2^{65}  &>& \frac{\sqrt{d_{66}}\cdot 9}{\pi[ \log d_{66}+\log 331]^{18}}\\
   &>& 7.2 \cdot 10^{20}.
\end{eqnarray*}
But $2^{65} \approx 3.7 \cdot 10^{19} < 7.2\cdot 10^{20}$, a contradiction.

\end{proof}
Other values for $A$ and $T$ can be chosen to obtain similar results.  Conrey and Iwaniec consider the case
$A=12$ and $\log T = (\log d)^{18}$ in Theorem \ref{CI Theorem} to get
\begin{equation*}
    L(1, \chi) \geq \frac{1}{(\log d)^{66}}.
\end{equation*}
With this lower bound a similar proof used for Theorem \ref{My Theorem} shows that there are no discriminants
$d$ with one class per genus for $d>d_{207} \approx 2.4 \cdot 10^{534}$.

We can also consider the case when $\psi = 1$ is the trivial class group character. Here the Riemann zeta
function appears as a factor of our Hecke L-function
\begin{equation*}
    L(s, \psi) = \zeta (s) L(s, \chi).
\end{equation*}
In this case Conrey and Iwaniec put a condition on the spacing of the zeros of $\zeta (s)$ (independent of $d$)
and use a more general version of Theorem \ref{CI Theorem} to get
\begin{equation*}
        L(1, \chi) \geq \frac{1}{(\log d)^{74}}
\end{equation*}
 \cite[p.309]{CI}. With this bound the method of Theorem
\ref{My Theorem} proves there are no discriminants with one class per genus for $d>d_{230} \approx 2.9 \cdot
10^{606}$.


\end{document}